\documentclass[12pt]{article}
\usepackage{latexsym}
\usepackage{amsmath,amssymb,amsthm}

\usepackage{amsmath,amsthm,latexsym,amssymb,pgf,makeidx,authblk,boxedminipage,tikz}
\usetikzlibrary{arrows,decorations.pathmorphing,backgrounds,positioning,fit,petri,calc}

\begin{document}

\newtheorem{defi}{Def}
\newtheorem{theo}{Theorem}
\newtheorem{coro}{Corollary}
\newtheorem{prop}{Proposition}
\newtheorem{rem}{Rem:}
\newtheorem{la}{Lemma}
\newtheorem{exa}{Example}
\newtheorem{conj}{Conjecture}
\newtheorem{problem}{Problem}

\title{On the Edge-Connectivity of the Square of a Graph}
\author{ Camino Balbuena\thanks{Departament de Matem\`{a}tica Aplicada III, 
Universitat Polit\`{e}cnica de Catalunya, Barcelona, Spain.} 
\and 
Peter Dankelmann\thanks{Department of Mathematics and Applied Mathematics, University of Johannesburg, Johannesburg, South Africa.} 
 } 
\maketitle

\begin{abstract}
Let $G$ be a connected graph. 
The edge-connectivity of $G$, denoted by $\lambda(G)$, is the minimum number of edges whose 
removal renders $G$ disconnected. Let $\delta(G)$ be the minimum degree of $G$. 
It is well-known that $\lambda(G) \leq \delta(G)$, and graphs for which equality holds are
said to be maximally edge-connected.
The square $G^2$ of $G$ is the graph with the same vertex set as $G$, in which two vertices
are adjacent if their distance is not more that $2$.

In this paper we present results on the edge-connectivity of the square of a graph.
We show that if the minimum degree of a connected graph $G$ of order $n$  is 
at least $\lfloor \frac{n+2}{4}\rfloor$, then $G^2$ is maximally edge-connected, and this
result is best possible. 
We also give lower bounds on $\lambda(G^2)$ for the case that  
$G^2$ is not maximally edge-connected: We prove that 
$\lambda(G^2) \geq \kappa(G)^2 + \kappa(G)$, where 
$\kappa(G)$ denotes the connectivity of $G$, i.e., the minimum number of 
vertices whose removal renders $G$ disconnected, and this bound is sharp. 
We further prove that $\lambda(G^2) \geq \frac{1}{2}\lambda(G)^{3/2} - \frac{1}{2} \lambda(G)$,
and we construct an infinite family of 
graphs to show that the exponent $3/2$ of $\lambda(G)$ in this bound is best possible. 
\end{abstract}


\section{Introduction}

The {\it edge-connectivity} $\lambda(G)$ of a graph $G$ is defined as the smallest number
of edges whose removal renders $G$ disconnected. It is one of the oldest
and best studied graph invariants. A classical result by Whitney \cite{Whi1932}
states that the edge-connectivity of a graph cannot exceed its {\it minimum degree},
defined as the smallest of the degrees of the vertices of $G$.  
Graphs whose edge-connectivity equals the minimum degree are called 
{\it maximally edge-connected} graphs. 
Much research on edge-connectivity has focused on 
finding sufficient conditions for a graph to be maximally edge-connected,
see for example the survey by Hellwig and Volkmann \cite{HelVol2008}.
Of particular relevance for us is the following classical 
result by Chartrand \cite{Cha1966}. 

\begin{theo}  \label{theo:Chartrand}
Let $G$ be a graph of order $n$ and minimum degree $\delta(G)$. If
\[ \delta(G) \geq \frac{n-1}{2}, \]
then $\lambda(G)=\delta(G)$. 
\end{theo}

Also results 
on the edge-connectivity of graphs derived from other graphs can be found in the literature, 
for example line graphs (see \cite{ChaSte1969, BauTin1979, HelRauVol2004, Men2001}), 
products of graphs (see \cite{UefVol2003}) and complements of graphs
(see, for example, \cite{AchAchCac1990,AlaMit1970, HelVol2008-2}).

The {\it square} $G^2$ of a graph $G$ is the graph on the
same vertex set as $G$, in which two vertices are adjacent if their
distance in $G$ is not more than two.  
Somewhat surprisingly, the edge-connectivity of the square of a 
graph appears not to have been studied. 
In this paper we aim to fill this gap in the literature. 

This paper is organised as follows. 
The terminology in this paper is defined in Scetion \ref{section:terminology-and-notation}. 
In Section \ref{section:sufficient condition for maximally edge-connected} 
we give a minimum degree condition on $G$ that guarantees 
that $G^2$ is maximally edge-connected;  
we prove that if a connected graph $G$ of order $n$ satisfies the
minimum degree condition $\delta(G) \geq \lfloor \frac{n+2}{4} \rfloor$,
then $\lambda(G^2)=\delta(G^2)$. We also construct a family of graphs to show
that this condition is best possible. 
In the following two sections we relate  connectivity-parameters of $G$ 
to the edge-connectivity of $G^2$.
In Section \ref{section:bounds in terms of kappa} we give a lower bound on 
$\lambda(G^2)$ in terms of $\kappa(G)$, the connectivity of $G$. We show that 
$\lambda(G^2) \geq \kappa(G)^2 + \kappa(G)$, unless $G^2$ is maximally edge-connected, 
and this bound is sharp.
In Section \ref{section:bounds in terms of lambda} we  give a lower bound on 
$\lambda(G)$ in terms of $\lambda(G)$. We prove that  
$\lambda(G^2) \geq \frac{1}{2}\lambda(G)^{3/2} - \frac{1}{2}\lambda(G)^{1/2}$, 
unless $G$ is maximally edge-connected. 
We construct an infinite family of graphs that show that the exponent 
$3/2$ of $\lambda(G)$ is best possible.

\section{Terminology and Notation}
\label{section:terminology-and-notation}

The notation we use is as follows. $G$ denotes a connected,
finite graph with vertex set $V(G)$ and edge set $E(G)$. 
The {\em order} of $G$, i.e., the number of vertices, is denoted by $n(G)$. 
An edge joining a vertex $u$ to a vertex $v$ is denoted by $uv$. 
If $v$ is a vertex of $G$
then we write $N_G(v)$ for the {\em neighbourhood} of $v$, i.e., the set 
of all vertices of $G$ adjacent to $v$, and $N_G[v]$ for the {\em closed 
neighbourhood} of $v$, i.e., the set $N_G(v) \cup \{v\}$. The {\em degree}
${\rm deg}_G(v)$ of $v$ is defined as $|N_G(v)|$, and the {\em minimum degree} of $G$,
denoted by $\delta(G)$, is $\min_{w\in V(G)}{\rm deg}_G(w)$. 
If there is no danger of confusion, then we often drop the argument or subscript $G$.

The {\em distance}  $d_G(u,v)$ 
between $u$ and $v$ is the minimum length of a $(u,v)$-path. 
An {\em edge-cut} of $G$ is a set of edges whose removal disconnects the graph, and 
a {\em minimum edge-cut} is an edge-cut of cardinality $\lambda(G)$.

If $A$ and $B$ are disjoint subsets of the vertex set of $G$, then $E_G(A,B)$
denotes the set of all edges of $G$ joining a vertex of $A$ to a vertex of $B$. 
The subgraph of $G$ {\em induced} by $A$,  denoted by $G[A]$, is the graph with vertex 
set $A$, in which two vertices are adjacent if and only if they are adjacent in $G$.

By $K_n$ we mean the complete graph on $n$ vertices.
The {\em sequential sum} of disjoint graphs $G_1, G_2, \ldots, G_k$, 
denoted by $G_1 + G_2 + \cdots + G_k$,  
is the graph obtained from their union by joining every vertex
in $G_i$ to every vertex in $G_{i+1}$ for $i=1,2,\ldots,k-1$. 

In the proofs below we use the following notation. 
Given a connected graph $G$, we denote the graph $G^2$ by $H$.
The set $S \subseteq E(H)$ is always a minimum edge-cut of
$H$. The two components of $H-S$ are denoted by $H_1$ and $H_2$,
and their vertex sets by $V_1$ and $V_2$, respectively. 
$A$ and $B$ are the sets of vertices in $H_1$ and $H_2$,
respectively, incident with edges in $S$. A vertex of $H_1$
or $H_2$ is an interior vertex if it is not in $A$ or $B$, respectively. 
We let $S' = S \cap E(G)$. Then $G-S'$ is disconnected since it
is an spanning subgraph of $H-S$. Let $G_1$ and $G_2$ be the
subgraphs of $G$ induced by $V_1$ and $V_2$. Then $G_1$ and
$G_2$ are unions of components of $G-S'$. Let $A'$ and $B'$ be the
set of vertices of $G_1$ and $G_2$, respectively, incident with
edges in $S'$. Note that $A' \subseteq A$ and $B' \subseteq B$.

\section{A sufficient condition for $G^2$ to be maximally edge-connected}
\label{section:sufficient condition for maximally edge-connected}

Theorem \ref{theo:Chartrand} shows that every graph $G$ of order $n$ with
$\delta(G) \geq \frac{n-1}{2}$ is maximally edge-connected. 
In this section we show that this condition can
be relaxed considerably to guarantee that $G^2$ is maximally
edge-connected. In our proof we make use of the following well-known lemma.

\begin{la}  {\rm (Hamidoune \cite{Ham1980})}
\label{la:hamidoune}
Let $G$ be a connected graph. If $S\subseteq E(G)$ is an edge-cut 
with $|S| < \delta(G)$,
then every component of $G-S$ has an interior vertex.
\end{la}

\begin{theo} \label{theo:mindegree-condition-for-G2-max-edge-conn}
Let $G$ be a connected graph of order $n$. If
\[ \delta(G) \geq \Big\lfloor \frac{n+2}{4} \Big\rfloor , \]
then $G^2$ is maximally edge-connected. 
\end{theo}

{\bf Proof:}  Suppose that $H$ is not maximally edge-connected. Let 
$H, S, H_1, H_2$, $S', A', B', V_1, V_2$, 
be as defined above. \\[1mm]
{\sc Claim 1:} $|S| \leq 2\delta(G)-2$. \\
Suppose to the contrary that 
$|S| \geq 2\delta(G)-1$. 
Applying the hypothesis of the theorem, we obtain that 
\[  |S| \geq 2\delta(G)-1
        \geq 2\big\lfloor \frac{n+2}{4}\big\rfloor-1
        \geq 2\frac{n-1}{4}-1
        = \frac{n-3}{2}. \]
Since $H$ is not maximally edge-connected, we thus have 
$\delta(H) \geq |S|+1 \geq \frac{n-1}{2}$. It follows by Theorem \ref{theo:Chartrand}
that $H$ is maximally edge-connected.        
This contradiction proves Claim 1. 
 \\[1mm]
{\sc Claim 2:} For all $u\in A'\cup B'$ we have 
$|N_G[u]\cap V_1|=1$ or $|N_G[u]\cap V_2|=1$ . \\
Clearly  $|N_G[u]\cap V_i|\geq 1$ for all $u\in A'\cup B'$ and $i\in \{1,2\}$. 
Suppose to the contrary that there exists $u\in A'\cup B'$ with 
 $|N_G[u]\cap V_1|\geq 2$ and $|N_G[u]\cap V_2|\geq 2$.
 Without loss of generality we may assume that $u\in A'$. 
Since $N_G[u]$ induces a complete graph of order ${\rm deg}_G(u)+1$ in $H$, we have that 
\begin{eqnarray}
|S \cap E(H[N_G[u]])| &\geq& |N_G[u] \cap V_1| |N_G[u] \cap V_2| \nonumber \\
     &= & |N_G[u] \cap V_1| ({\rm deg}_G(u)+1-|N_G[u] \cap V_1|)  \nonumber \\
      &\geq&  2 ({\rm deg}_G(u)-1) \nonumber \\
      & \geq & 2\delta(G) -2.  \label{eq:lower-bound-on-S}
\end{eqnarray}       
From \eqref{eq:lower-bound-on-S} and Claim 1 we obtain that 
$|S \cap E(H[N_G[u]])|=|S|=2\delta(G)-2$,  
and so 
\begin{equation} \label{eq:S-contained-in-N[u]}
S \subseteq E(H[N_G[u]]). 
\end{equation}
This implies in particular that $B' \subseteq N_G(u)$.
By Lemma \ref{la:hamidoune}, $H_2$ contains an internal vertex $v_2$. Let 
$Q$ be a shortest path in $G$ from $B'$ to $v_2$. Let $w \in B'$ be the initial vertex
of $Q$ and  $x$ its successor. Then $ux \in E(H)$. Since clearly, $x \in V_2$,
we have $ux \in S$.  But $x \notin N_G[u]$, and so $ux \notin E(H[N_G[u]])$. 
This contradiction to \eqref{eq:S-contained-in-N[u]} proves Claim 2.  \\[1mm]
{\sc Claim 3:} There exist $u\in A', v \in B'$ with
$|N_G[u] \cap V_2|=|N_G[v]\cap V_1|=1$. \\
By Lemma \ref{la:hamidoune} there exists an internal vertex
$v_1$ of $H_1$, i.e., a vertex not in $A$ and thus not in $A'$. 
Let $P$ be a shortest path in $G$ from $v_1$ to $A'$, let
$u$ be the terminal vertex of $P$ and $y$ its predecessor. 
Then $y\in V_1$, and so $\{u,y\} \subseteq N_G[u] \cap V_1$. By Claim 1 we have
$|N_G[u] \cap V_2| =1$, as desired. Similarly we 
show the existence of a vertex $v\in B'$ with
$|N_G[v]\cap V_1|=1$. Claim 3 follows. 

We are now ready to complete the proof. 
Let $u$ and $v$ be as in Claim 3, and let 
$u^*$ and $v^*$ be the unique neighbour of 
$u$ and $v$ in $V_2$ and $V_1$, respectively, in $G$.
Let $E_u \subseteq E(H)$ be the set of all edges of the
form $xu^*$, where $x\in N_G[u]\cap V_1$, and similarly
let $E_v \subseteq E(H)$ be the set of all edges of the
form $xv^*$, where $x\in N_G[v]\cap V_2$.
Clearly,
\[ E_u \cup E_v \subseteq S. \]
Clearly, $|E_u| = {\rm deg}_G(u) \geq \delta(G)$, 
and similarly $|E_v| \geq \delta(G)$. Since all edges
in $E_u$ are incident with $u^*$, and all edges in
$E_v$ are incident with $v^*$, we conclude that
$E_u \cap E_v$ is either empty or contains the edge $u^*v^*$
but no further edge. Hence
\[ |S| \geq |E_u| + |E_v| - |E_u\cap E_v| 
        \geq 2\delta(G)-1. \]
This contradiction to Claim 2 concludes the proof 
of the theorem. 
\hfill $\Box$ \\

The condition $\delta(G)\geq \lfloor \frac{n+2}{4} \rfloor$ in Theorem 
\ref{theo:mindegree-condition-for-G2-max-edge-conn} is best possible for
all $n$ with $n\geq 10$. For given $n$ let $G_n$ be the sequential sum
\[ K_{n/4} + K_{(n-4)/4} + K_1 + K_1 +  K_{(n-8)/4} + K_{(n+4)/4} \quad 
                   \textrm{if $n\equiv 0 \pmod 4$,}  \]
\[ K_{(n+3)/4} + K_{(n-5)/4} + K_1 + K_1 +  K_{(n-9)/4} + K_{(n+3)/4} \quad 
                   \textrm{if $n\equiv 1 \pmod 4$,}  \]
\[ K_{(n+2)/4} + K_{(n-6)/4} + K_1 + K_1 +  K_{(n-6)/4} + K_{(n+2)/4} \quad 
                     \textrm{if $n\equiv 2 \pmod 4$,}  \]
\[ K_{(n+1)/4} + K_{(n-7)/4} + K_1 + K_1 +  K_{(n-7)/4} + K_{(n+5)/4} \quad 
                    \textrm{if $n\equiv 3 \pmod 4$.}  \]
It is easy to verify that $\delta(G_n)=\lfloor\frac{n+2}{4}\rfloor -1$, 
but $\lambda(G_n^2)=\delta(G_n^2)-1$. \\

\section{A lower bound on $\lambda(G^2)$ in terms of $\kappa(G)$}
\label{section:bounds in terms of kappa}

In this section we explore the relationship between the connectivity of a 
graph and the edge-connectivity of its square. For graphs whose square is not 
maximally edge-connected, we give a lower bound on the  edge-connectivity 
of $G^2$ in terms of the minimum degree and the connectivity of $G$.
As a corollary we obtain a sharp lower bound on the  edge-connectivity 
of $G^2$ in terms of the connectivity of $G$ only.

\begin{theo} \label{theo:kappa-delta}
Let $G$ be a connected graph with $\delta(G)\geq 2$. Then
$G^2$ is maximally edge-connected, or 
\[ \lambda(G^2) \geq \kappa(G)(\delta(G)+1). \]
\end{theo}

{\bf Proof:}
Let $H$, $S$, $H_1$, $H_2$ and $S'$ be as defined above. 
We may assume that $H$ is not maximally edge-connected, otherwise
the theorem holds. 
By Lemma \ref{la:hamidoune} there exist interior vertices $v_1$ of $H_1$ and
$v_2$ of $H_2$. 

Let $\kappa := \kappa(G)$ and $\delta := \delta(G)$. By Menger's
Theorem, there exist $\kappa$ internally vertex-disjoint 
$(v_1,v_2)$-paths $P_1, P_2, \ldots,P_{\kappa}$ in $G$. 
We may assume that each $P_i$ is an induced path in $G$. 
Each 
$P_i$ uses at least one edge of $S'$. Let $a_ib_i$ be the 
last edge of $P_i$ in $S'$.

We now show that 
\begin{equation}  \label{eq:kappa-1}
|S| \geq \kappa + \sum_{i=1}^{\kappa} {\rm deg}_G(a_i). 
\end{equation}
Let $M_1 := \{a_1, a_2,\ldots,a_{\kappa}\}$. 
Consider the edges of $G$ incident with a vertex in $M_1$. 
First let $e$ be an edge incident with exactly one vertex in $M_1$, so $e=a_ix$ for 
some $i\in \{1,2,\ldots,\kappa\}$ and $x \in V(G) - M_1$.  
From $e$ we obtain an edge $e_S \in S$  as follows.
If $x \in V_2$, then let $e_S=e$. 
If $x \in V_1-M_1$, then let $e_S=b_ix$. 
Now let $e$ be an edge incident with two vertices in $M_1$, so 
$e=a_ia_j$ for some $i,j \in \{1,2,\ldots,\kappa\}$. We 
obtain two edges from $e$, the edges $a_iy_j$ and $a_jy_i$, where 
$y_j$ is the first vertex of $P_j$ after $a_j$ that is not in $N_G(a_i)$,
and $y_i$ is the first vertex of $P_i$ after $a_i$ that is not in $N_G(a_j)$. 
Clearly, $y_i, y_j \in V_2$, so $a_iy_j, a_jy_i \in S$. 
Note that $y_i,y_j \neq v_2$ since otherwise $v_2$ would not be an internal vertex of $H$, 
and so $y_i \neq y_j$.

Let $E_1$ be the set of edges of $S$ thus obtained from the edges incident with vertices of 
$M_1$. Then each edge in $E_1$ is either
(i) of the form $xb_i$, where $x \in V_1-M_1$, in which case it is obtained from the edge $a_ix$, or (ii) of the form $a_ix$, where $x \in V_2$ and $a_ix$ is an edge of $G$, in which case it 
is obtained from $a_ix$, or
(iii) of the form $a_ix$, where $x \in V_2 \cap V(P_j)$ and $a_ix$ is not an edge of $G$. 
It follows that we do not obtain the same edge from two distinct edges of $G$ incident 
with a vertex in $M_1$. Therefore, 
\[   |E_1| = \sum_{i=1}^{\kappa} {\rm deg}_G(a_i). \]
For $i=1,2,\ldots,\kappa$ let $w_i$ be the vertex following $b_i$ on $P_i$. 
Then $a_iw_i \in S$, but it follows from the above that $a_iw_i \notin E_1$, so 
$S$ contains at least $\kappa$ edges in addition to the edges in $E_1$. 
Hence \eqref{eq:kappa-1} follows.

Since each vertex of $M_1$ has degree at least $\delta$ in $G$, \eqref{eq:kappa-1} 
yields that
\begin{equation} \label{eq:kappa-2}
|S| \geq |E_1|+ \kappa \geq \kappa + \kappa\delta. 
\end{equation}                  
Since $|S| = \lambda(G^2)$, the theorem follows. \hfill $\Box$ \\

We do not know if the bound in Theorem \ref{theo:kappa-delta} is sharp.
However, using the inequality $\delta(G) \geq \kappa(G)$, we obtain the 
following sharp lower bound on the edge-connectivity of $G^2$ in terms of
the connectivity of $G$ as a corollary to Theorem \ref{theo:kappa-delta}.

\begin{coro}  \label{coro:bound-in-terms-of-kappa}
Let $G$ be a $\kappa$-connected graph with $\kappa \geq 2$. Then $G^2$ is 
maximally edge-connected or
\[ \lambda(G^2) \geq \kappa(\kappa+1). \]
If $\kappa$ is even, then the bound is sharp. 
\end{coro}

{\bf Proof:}
The bound in Corollary \ref{coro:bound-in-terms-of-kappa} follows from 
Theorem \ref{theo:kappa-delta} in conjunction with the inequality 
$\kappa(G) \leq \delta(G)$. 

\begin{figure}[h]
  \begin{center}
\begin{tikzpicture}
  [scale=0.6,inner sep=1mm, 
   vertex/.style={circle,thick,draw}, 
   thickedge/.style={line width=2pt}] 
    \node[vertex] (a1) at (4,0) [fill=white] {};
    \node[vertex] (a2) at (7,0) [fill=white] {};
    \node[vertex] (a3) at (10,0) [fill=white] {};
    \node[vertex] (a4) at (13,0) [fill=white] {};       
    \node[vertex] (b1) at (4,2) [fill=white] {};
    \node[vertex] (b2) at (7,2) [fill=white] {};
    \node[vertex] (b3) at (10,2) [fill=white] {};
    \node[vertex] (b4) at (13,2) [fill=white] {};      
    \node[vertex] (c1) at (4,4) [fill=white] {};
    \node[vertex] (c2) at (7,4) [fill=white] {};
    \node[vertex] (c3) at (10,4) [fill=white] {};
    \node[vertex] (c4) at (13,4) [fill=white] {};   
    \node[vertex] (d1) at (4,6) [fill=white] {};
    \node[vertex] (d2) at (7,6) [fill=white] {};
    \node[vertex] (d3) at (10,6) [fill=white] {};
    \node[vertex] (d4) at (13,6) [fill=white] {};

    \draw[black, thick] (a1)--(a2)--(a3)--(a4);
    \draw[black, thick] (b1)--(b2)--(b3)--(b4);
    \draw[black, thick] (c1)--(c2)--(c3)--(c4);        
    \draw[black, thick] (d1)--(d2)--(d3)--(d4);

   \draw[black, dotted, thick] (a1)--(b1)  (c1)--(d1);
   \draw[black, dotted, thick] (a2)--(b2)  (c2)--(d2);
   \draw[black, dotted, thick] (a3)--(b3)  (c3)--(d3);
   \draw[black, dotted, thick] (a4)--(b4)  (c4)--(d4);

   \draw[black, thick, rounded corners] (0,-0.4) rectangle (2,6.4);      
   \node at (1,3) {$K_{2\kappa^2}$};   
   \draw[gray,  rounded corners] (3.5,-0.4) rectangle (4.5,6.4);    
   \draw[gray,  rounded corners] (6.5,-0.4) rectangle (7.5,6.4);    
   \draw[gray,  rounded corners] (9.5,-0.4) rectangle (10.5,6.4);    
   \draw[gray,  rounded corners] (12.5,-0.4) rectangle (13.5,6.4);                
   \draw[black, thick, rounded corners] (15,-0.4) rectangle (17,6.4);      
   \node at (16,3) {$K_{2\kappa^2}$};   
   
   \node at (4,-1) {$K_{\kappa}-M$};      
   \node at (7,-1) {$K_{\kappa}-M$};      
   \node at (10,-1) {$K_{\kappa}-M$};   
   \node at (13,-1) {$K_{\kappa}-M$};   
          
   \draw[black, thick] (2,2.9)--(3.5,2.9)  (2,3.1)--(3.5,3.1);
   \draw[black, thick] (13.5,2.9)--(15,2.9)  (13.5,3.1)--(15,3.1);

\end{tikzpicture}
\end{center}
\caption{The graph $G(\kappa)$ for $\kappa=4$. Dotted lines indicate edges not present.}
\label{fig:graph-G4}
\end{figure}
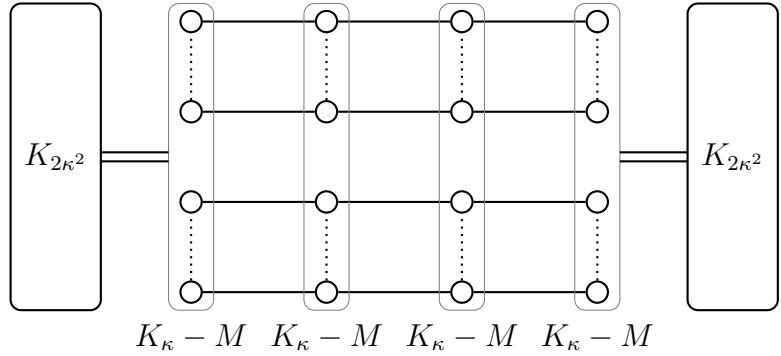

It remains to prove the sharpness of the bound. 
For even $\kappa \in \mathbb{N}$ consider the graph $K_{\kappa}-M$, where $K_{\kappa}$
is a complete graph with vertex set $\{u_1, u_2,\ldots, u_{\kappa}\}$, and 
$M$ is the perfect matching $\{u_1u_2, u_3u_4,\ldots, u_{\kappa-1}u_{\kappa}\}$ of $K_{\kappa}$. 
For $i=1,2,3,4$ let $G_i$ be a copy of $K_{\kappa}-M$, and denote the vertex of $G_i$ 
corresponding to a vertex $u_j$ of $K_{\kappa}-M$ by $u_j^{(i)}$.
Let $G_0$ and $G_5$ be disjoint copies of the complete graph $K_{2\kappa^2}$. 
Then we obtain the graph $G(\kappa)$ from joining every vertex of $G_0$ to every vertex
of $G_1$, and every vertex of $G_4$ to every vertex of $G_5$ by an edge, and by 
adding the edges $v_i^{(1)}v_i^{(2)}, v_i^{(2)}v_i^{(3)}$ and $v_i^{(3)}v_i^{(4)}$ for 
$i=1,2,\ldots,\kappa$. 
The graph $G(\kappa)$ for $\kappa=4$ is sketched in Figure \ref{fig:graph-G4}.

It is easy to verify that $G(\kappa)$ is $\kappa$-connected. Denote the vertex set of 
$G_i$ by $V^i$. Clearly, removing the edges in 
$E_{G^2}(V^0 \cup V^1 \cup V^2, V^3\cup V^4 \cup V^5)$ disconnects the graph $G^2$.
Now $|E_{G^2}(V^2,V^3)|=\kappa(\kappa-1)$, 
$|E_{G^2}(V^1,V^3)|=|E_{G^2}(V^2,V^4)|=\kappa$, and so 
\begin{eqnarray*} 
\lambda(G^2) & \leq & |E_{G^2}(V^0 \cup V^1 \cup V^2, V^3\cup V^4 \cup V^5)| \\
      &=& |E_{G^2}(V^2,V^3)|+  |E_{G^2}(V^1,V^3)|+ |E_{G^2}(V^2,V^4)|   \\
      & = & \kappa(\kappa+1). 
\end{eqnarray*}      
Since every vertex of $G^2$ is adjacent to all vertices of $V^0$ or all vertices of $V^5$,
except possibly itself, we have $\delta(G^2) \geq 2\kappa^2 > \lambda(G^2)$, so 
$G^2$ is not maximally edge-connected. Hence Corollary \ref{coro:bound-in-terms-of-kappa} 
yields that $\lambda(G^2) \geq \kappa(\kappa+1)$, so $G^2$ attains equality in 
Corollary \ref{coro:bound-in-terms-of-kappa}.
\hfill $\Box$  \\

\begin{figure}[h]
  \begin{center}
\begin{tikzpicture}
  [scale=0.6,inner sep=1mm, 
   vertex/.style={circle,thick,draw}, 
   thickedge/.style={line width=2pt}] 
    \node[vertex] (a1) at (4,0) [fill=white] {};
    \node[vertex] (a2) at (7,0) [fill=white] {};
    \node[vertex] (a3) at (10,0) [fill=white] {};
    \node[vertex] (a4) at (13,0) [fill=white] {};       
    \node[vertex] (b1) at (4,2) [fill=white] {};
    \node[vertex] (b2) at (7,2) [fill=white] {};
    \node[vertex] (b3) at (10,2) [fill=white] {};
    \node[vertex] (b4) at (13,2) [fill=white] {};      
    \node[vertex] (c1) at (4,4) [fill=white] {};
    \node[vertex] (c2) at (7,4) [fill=white] {};
    \node[vertex] (c3) at (10,4) [fill=white] {};
    \node[vertex] (c4) at (13,4) [fill=white] {};   
    \node[vertex] (d1) at (4,6) [fill=white] {};
    \node[vertex] (d2) at (7,6) [fill=white] {};
    \node[vertex] (d3) at (10,6) [fill=white] {};
    \node[vertex] (d4) at (13,6) [fill=white] {};      
    \node[vertex] (e1) at (4,8) [fill=white] {};
    \node[vertex] (e2) at (7,8) [fill=white] {};
    \node[vertex] (e3) at (10,8) [fill=white] {};
    \node[vertex] (e4) at (13,8) [fill=white] {};

    \draw[black, thick] (a1)--(a2)--(a3)--(a4);
    \draw[black, thick] (b1)--(b2)--(b3)--(b4);
    \draw[black, thick] (c1)--(c2)--(c3)--(c4);        
    \draw[black, thick] (d1)--(d2)--(d3)--(d4);
    \draw[black, thick] (e1)--(e2)--(e3)--(e4);    
    \draw[black, thick] (a1)--(b2)  (b3)--(a4);

   \draw[black, dotted, thick] (a1)--(b1)--(c1)  (d1)--(e1);
   \draw[black, dotted, thick] (a2)--(b2)--(c2)  (d2)--(e2);
   \draw[black, dotted, thick] (a3)--(b3)--(c3)  (d3)--(e3);
   \draw[black, dotted, thick] (a4)--(b4)--(c4)  (d4)--(e4);

   \draw[black, thick, rounded corners] (0,-0.4) rectangle (2,8.4);      
   \node at (1,4) {$K_{2\kappa^2}$};   
   \draw[gray,  rounded corners] (3.5,-0.4) rectangle (4.5,8.4);    
   \draw[gray,  rounded corners] (6.5,-0.4) rectangle (7.5,8.4);    
   \draw[gray,  rounded corners] (9.5,-0.4) rectangle (10.5,8.4);    
   \draw[gray,  rounded corners] (12.5,-0.4) rectangle (13.5,8.4);                
   \draw[black, thick, rounded corners] (15,-0.4) rectangle (17,8.4);      
   \node at (16,4) {$K_{2\kappa^2}$};   
   
   \node at (4,-1) {$K_{\kappa}-M$};      
   \node at (7,-1) {$K_{\kappa}-M$};      
   \node at (10,-1) {$K_{\kappa}-M$};   
   \node at (13,-1) {$K_{\kappa}-M$};   
          
   \draw[black, thick] (2,3.9)--(3.5,3.9)  (2,4.1)--(3.5,4.1);
   \draw[black, thick] (13.5,3.9)--(15,3.9)  (13.5,4.1)--(15,4.1);

\end{tikzpicture}
\end{center}
\caption{The graph $G(\kappa)$ for $\kappa=5$. Dotted lines indicate edges not present.}
\label{fig:graph-G5}
\end{figure}
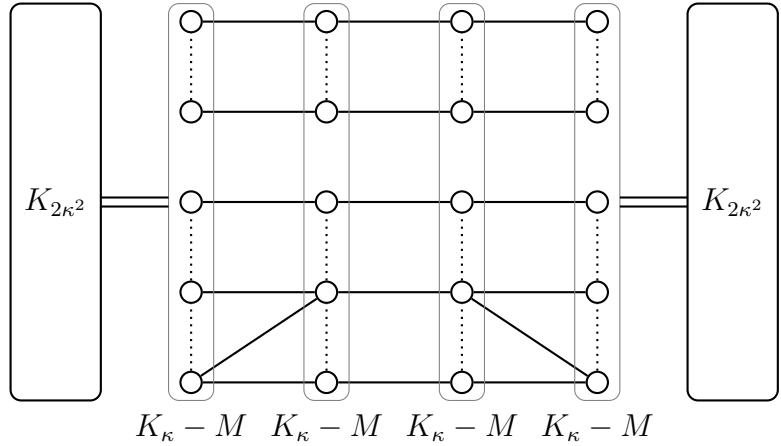

We note that for odd $\kappa$, there exist graphs for which the edge-connectivity 
of the square differs from the bound in Corollary \ref{coro:bound-in-terms-of-kappa} by $1$. 
To see that consider, for odd $\kappa \in \mathbb{N}$, the graph $K_{\kappa}-M_0$, where $K_{\kappa}$
is a complete graph with vertex set $\{u_1, u_2,\ldots, u_{\kappa}\}$, and 
$M$ is the minimum edge-cover $\{u_1u_2, u_3u_4,\ldots, u_{\kappa-2}u_{\kappa-1}, 
u_{\kappa-1}u_{\kappa}$ of $K_{\kappa}$. 
For $i=1,2,3,4$ let $G_i$ be a copy of $K_{\kappa}-M$, and denote the vertex of $G_i$ 
corresponding to a vertex $u_j$ of $K_{\kappa}-M$ by $u_j^{(i)}$.
Let $G_0$ and $G_5$ be disjoint copies of the complete graph $K_{2\kappa^2}$. 
Then we obtain the graph $G(\kappa)$ from joining every vertex of $G_0$ to every vertex
of $G_1$, and every vertex of $G_4$ to every vertex of $G_5$ by an edge, and by 
adding the edges $v_i^{(1)}v_i^{(2)}, v_i^{(2)}v_i^{(3)}$ and $v_i^{(3)}v_i^{(4)}$ for 
$i=1,2,\ldots,\kappa$. Finally we add the edges $u_{\kappa-1}^{(2)} u_{\kappa}^{(1)}$ 
and $u_{\kappa-1}^{(3)} u_{\kappa}^{(4)}$. 
The graph $G(\kappa)$ for $\kappa=5$ is sketched in Figure \ref{fig:graph-G5}.
It is easy to verify that $G(\kappa)$ is $\kappa$-connected. Similar 
considerations as for the case $\kappa$ even show that 
\[   \lambda(G(\kappa)^2) \leq \kappa(\kappa+1)+1. \]

\begin{coro}
Let $G$ be a connected graph with minimum degree $\delta$. If $G$ is not complete, then
\[ \lambda(G^2) \geq \delta+1, \]
and this bound is sharp. 
\end{coro}

{\bf Proof:}
If $G$ is not complete, then clearly $\delta(G^2)\geq \delta+1$. Thus, if $G^2$ is 
maximally edge-connected, then the statement of the corollary holds. If $G^2$ is not
maximally edge-connected, then $\kappa(G)\geq 1$, and the corollary follows from 
Theorem \ref{theo:kappa-delta}.

To see that the bound is sharp consider for given $\delta \in \mathbb{N}$ the 
graph $G= K_1 + K_{\delta} + K_{1} + K_{\delta}$.  Clearly $\delta(G) = \delta$ 
and $\lambda(G^2)= \delta(G^2) = \delta+1$.
\hfill $\Box$

\section{A lower bound on $\lambda(G^2)$ in terms of $\lambda(G)$}
\label{section:bounds in terms of lambda}

In this section we explore the relationship between the edge-connectivity of $G$ and
the edge-connectivity of $G^2$.  We show that $\lambda(G)$ is bounded from below by 
a term $O(\lambda(G)^{3/2})$, unless $G^2$ is maximally edge-connected.
We also construct an infinite family of graphs that show that 
the exponent $3/2$ is best possible. 

\begin{theo}   \label{theo:lambda}
Let $G$ be a $\lambda$-edge-connected graph. Then $G^2$ is maximally edge-connected or
\[ \lambda(G^2) \geq \frac{1}{2}\lambda^{3/2} - \frac{1}{2}\lambda. \] 
\end{theo}

{\bf Proof:}
Let $H, S, A, B, S', A',B'$ be as defined above. For $v \in A' \cup B'$
we also define ${\rm inc}(v,S')$ to be the number of edges in $S'$ 
incident with $v$. 

We may assume that $H$ is not maximally edge-connected, otherwise the theorem holds. 
\\[1mm]
{\sc Case 1:} There exists a vertex $v \in A' \cup B'$ with 
$\sqrt{\lambda} \leq {\rm inc}(v,S') \leq \lambda - \lambda^{1/2}$. \\
We may assume that $v \in A'$. 
If $u \in N_G(v) \cap V_1$ and $w \in N_G(v) \cap V_2$, then the edge $uw$ of $H$ is in $S$. 
Hence
\begin{eqnarray} 
|S| & \geq & |N_G(v) \cap V_1| \cdot |N_G(v) \cap V_2| \nonumber \\
    & = & \big( {\rm deg}_G(v) - {\rm inc}(v,S') \big) {\rm inc}(v,S')  \nonumber \\
    & \geq &  \big( \lambda - {\rm inc}(v,S') \big) {\rm inc}(v,S'). \nonumber
\end{eqnarray}
For a real $x$ with $\lambda \leq x \leq \lambda - \lambda^{1/2}$, the function
$(\lambda-x)x$ is minimised if $x= \sqrt{\lambda}$ or $x = \lambda - \lambda^{1/2}$. 
Substituting either of these values yields
$\lambda(H) = |S| \geq \lambda^{3/2} - \lambda$, as desired.\\[1mm]
{\sc Case 2:} For every $v \in A'\cup B'$ we have ${\rm inc}(v,S') < \sqrt{\lambda}$
or ${\rm inc}(v,S') > \lambda - \sqrt{\lambda}$. \\
We denote by $A^+$  and $A^-$ the set of all $v\in A'$ with 
${\rm inc}(v,S') > \lambda - \sqrt{\lambda}$ and
${\rm inc}(v,S') < \sqrt{\lambda}$, respectively. The sets $B^+$ and $B^-$ are defined analogously. 
Then $A^+ \cup A^- = A'$, $B^+ \cup B^- = B'$. 

Define the sets $W_1$ and $W_2$ by $W_1 = (V_1-A^+) \cup B^+$ and 
$W_2 = (V_2  - B^+) \cup A^+$ and consider the edge-cut $T := E_G(W_1,W_2)$ of $G$. Then
\begin{equation} \label{eq:lambda-define-T}
T = E_G(A^+, B^+) \cup E_G(A^-, B^-) \cup E_G(A^+, V_1-A^+) \cup E_G(B^+, V_2 - B^+). 
\end{equation} 
Since $G$ is $\lambda$-edge-connected, we have $|T| \geq \lambda$, and so at least one of
the terms
$|E_G(A^+, B^+)|$,  $|E_G(A^-, B^-)|$,  $|E_G(A^+, V_1-A^+)|$ and $|E_G(B^+, V_2 - B^+)|$ 
is greater or equal to $\frac{1}{4} \lambda$. We consider the following (not mutually 
exclusive) subcases. \\[1mm]
{\sc Case 2a:} $|E_G(A^+, B^+)| \geq \frac{1}{4}\lambda$. \\[1mm]
Since $|E_G(A^+, B^+)| \leq |A^+| \cdot |B^+|$, we have 
$|A^+| \cdot |B^+| \geq \frac{1}{4}\lambda$ and so $|A^+| + |B^+| \geq \sqrt{\lambda}$. 
We may assume that $|A^+| \geq |B^+|$, and thus $|A^+| \geq \frac{1}{2} \sqrt{\lambda}$.
Since each vertex of $A^+$ is incident with at least $\lambda-\sqrt{\lambda}$ edges of $S'$,
we have
\begin{equation} \label{eq:lambda-case2a} 
|S'| \geq \sum_{v \in A^+} {\rm inc}(v,S') 
     \geq \frac{1}{2} \sqrt{\lambda} \big( \lambda- \sqrt{\lambda} \big)
     = \frac{1}{2} \lambda^{3/2} - \frac{1}{2}\lambda. 
\end{equation}     
Since $\lambda(H) = |S| \geq |S'|$, the theorem holds in this case. \\[1mm]  
{\sc Case 2b:} $|E_G(A^-, B^-)| \geq \frac{1}{4}\lambda$. \\[1mm]    
Let $A^*$ ($B^*$) be the set of vertices in $A^-$ that are adjacent in $G$ to a vertex in 
$B^-$ ($A^-$). The defining condition of this case yields that 
$|E_G(A^*, B^*)| \geq \frac{1}{4} \lambda$. We may assume that $|A^*| \geq |B^*|$, 
and the same arguments as in Case 2a show that
$|A^*| \geq \frac{1}{2} \sqrt{\lambda}$.  
Each vertex $v \in A^*$ has a neighbour $w$ in $B^-$, and it follows from the definition
of $B^-$ that $w$ has at least $\lambda - \sqrt{\lambda}$ neighbours in $V_2$, and
$v$ is joined to each of these by an edge in $H$, which is in $S$. Hence  
\[ |S| \geq \sum_{v \in A^*} {\rm inc}(v, S) 
             \geq \frac{1}{2}\sqrt{\lambda} \big( \lambda - \sqrt{\lambda} \big)
               = \frac{1}{2} \lambda^{3/2} - \frac{1}{2} \lambda, \]
and the theorem holds also in this case. \\[1mm]   
{\sc Case 2c:} $|E_G(A^+, V_1-A^+)| \geq \frac{1}{4} \lambda$ or 
 $|E_G(B^+, V_2 - B^+)| \geq \frac{1}{4} \lambda$. \\
We may assume that $|E_G(A^+, V_1-A^+)| \geq \frac{1}{4} \lambda$. Moreover, 
we may assume that $|A^+| < \frac{1}{2} \sqrt{\lambda}$ since otherwise, if 
$|A^+| \geq \frac{1}{2} \sqrt{\lambda}$, the same reasoning as in \eqref{eq:lambda-case2a}
shows that the theorem holds. 

Now let $v_0 \in A^+$ be a vertex with the most neighbours in $V_1-A^+$. Then 
\[ |N_G(v_0, V_1-A^+)| \geq \frac{ |E_G(A^+, V_1-A^+)| }{|A^+|}
                  \geq \frac{ \frac{1}{4} \lambda }{\frac{1}{2} \sqrt{\lambda} }
                  = \frac{1}{2} \sqrt{\lambda}. \]
Since every vertex in $N_G(v_0) \cap V_1$ is joined to every vertex in $N_G(v_0) \cap V_2$,
by an edge in $S$, we have
\[ |S| \geq |N_G(v_0, V_1-A^+)| \cdot {\rm inc}(v_0,S')
             \geq \frac{1}{2} \sqrt{\lambda} \cdot \big( \lambda - \sqrt{\lambda} \big)
             = \frac{1}{2}\lambda^{3/2} - \frac{1}{2} \lambda. \]
In all cases, the theorem holds. \hfill $\Box$  \\

The exponent $3/2$ of $\lambda$ in the bound in Theorem \ref{theo:lambda} is 
best possible. To see this consider for given square integer $\lambda \geq 4$ the
graph $G_{\lambda}$, defined by
\[ G_{\lambda} = K_{\lambda^2} + K_{\lambda- \sqrt{\lambda}} + K_{\sqrt{\lambda}} + 
    K_{\sqrt{\lambda}} + K_{\lambda- \sqrt{\lambda}} + K_{\lambda^2}. \]
It is easy to verify that $G_{\lambda}$ is $\lambda$-edge-connected, that
$\delta(G_{\lambda}^2)= \lambda^2+ \lambda-1$, and that 
$\lambda(G_{\lambda}^2)=2\lambda^{3/2}-\lambda$.

\section{Conclusion and open problems}

In Section \ref{section:sufficient condition for maximally edge-connected}, we presented
a degree condition on a graph $G$ that guarantees that $G^2$ is maximally edge-connected. 
Our condition in Theorem \ref{theo:mindegree-condition-for-G2-max-edge-conn} is a 
relaxation of the condition $\delta(G) \geq \frac{n-1}{2}$ in Theorem \ref{theo:Chartrand}. 
A is natural to ask if other sufficient conditions for maximally
edge-connected graphs can be relaxed in similar ways to guarantee that the square of a
graph is maximally edge-connected. Possible candidates are  
minimum degree conditions for graphs not containing certain subgraphs (see, for
example, \cite{Dan2019}),
degree sequence conditions (see, for example, \cite{DanVol1997}), 
conditions in terms of the inverse degree, i.e., the
sum of the reciprocals of the vertex degrees, (see \cite{DanHelVol2009}),
conditions in terms of size \cite{VolHon2017}, or 
conditions in terms of girth and diameter (see, for example \cite{SonNakIma1985}).

It is natural to expect that
the lower bounds on the edge-connectivity of $G^2$ in terms of the connectivity of $G$ in 
Section \ref{section:bounds in terms of kappa} and 
in terms of the edge-connectivity of $G$ in Section \ref{section:bounds in terms of lambda}
can be generalised to higher powers of graphs. For $k\in \mathbb{N}$, the $k$-th power of 
$G$ is the graph on the same vertex set as $G$, where two vertices are adjacent if their 
distance is not more than $k$. It would also be interesting to determine
the smallest minimum degree of a graph $G$ of given order  that guarantees that its $k$-th power
is maximally edge-connected. We pose these questions as open problems.


\begin{thebibliography}{99}
\bibitem{AchAchCac1990} N.\ Achuthan, N.R.\ Achuthan, L.\ Caccetta, 
   On the Nordhaus–Gaddum class problems.
   Australasian J.\ Combin.\ {\bf 2} (1990), 5-27.
\bibitem{AlaMit1970} Y.\ Alavi, J.\ Mitchem, 
   Connectivity and line connectivity of complementary graphs, 
   in: M. Capobianco, J.B. Frenchen, M. Krolik (Eds.), Recent Trends
   in Graph Theory (1970), 1-3.   
\bibitem{BauTin1979} D.\ Bauer, R.\ Tindell, 
   Graphs with prescribed connectivity and line graph connectivity. 
   J.\ Graph Theory {\bf 3} (1979), 393-395.   
\bibitem{Cha1966} G.\ Chartrand, 
   A graph theoretic approach to a communications problem. 
   SIAM J.\ Appl.\ Math.\ {\bf 14} (1966), 778-781.
\bibitem{ChaSte1969} G.\ Chartrand, M.J.\ Stewart, 
   The connectivity of line-graphs. 
   Math.\ Ann.\ {\bf 182} (1969), 170-174.  
\bibitem{DanVol1997} P.\ Dankelmann,  L.\ Volkmann,    
   Degree sequence conditions for maximally edge-connected graphs and digraphs.
   J.\ Graph Theory {\bf 26} no.\ 1 (1997), 27-34.
\bibitem{DanHelVol2009} P.\ Dankelmann, A.\ Hellwig, L.\ Volkmann,
   Inverse degree and edge-connectivity. 
   Discrete Math.\ {\bf 309}(9), (2009), 670-673.   
\bibitem{Dan2019} P.\ Dankelmann, 
   On the edge-connectivity of $C_4$-free graphs. 
   Commun.\ Comb.\ Optim.\ {\bf 4}(2) (2019), 141-150.          
\bibitem{Ham1980} Y.O.\ Hamidoune, 
   A property of $a$-fragments of a digraph.
   Discrete Math.\ {\bf 3} no.\ 1 (1980), 105-106. 
\bibitem{HelRauVol2004} A.\ Hellwig, D.\ Rautenbach, L.\ Volkmann, 
   Note on the connectivity of line graphs. 
   Inform.\ Process.\ Lett.\ {\bf 91} (2004), 7-10.   
\bibitem{HelVol2008} A.\ Hellwig, L.\ Volkmann, 
   Maximally edge-connected and vertex-connected graphs and digraphs - a survey.
   Discrete Math.\ {\bf 308} (2008), 3265-3296.   
\bibitem{HelVol2008-2} A.\ Hellwig, L.\ Volkmann,    
   The connectivity of a graph and its complement. 
   Discrete Appl.\ Math. {\bf 156} (2008), 3325-3328.
\bibitem{Men2001} J.\ Meng, 
   Superconnectivity and super edge-connectivity of line graphs. 
   Graph Theory Notes N.\ Y.\ {\bf 40} (2001), 12-14.   
\bibitem{SonNakIma1985} T.\ Soneoka, H.\ Nakada, M.\ Imase, 
   Sufficient conditions for dense graphs to be maximally connected. 
   Discrete Math.\ {\bf 63}(1) (1987), 53-66.
\bibitem{UefVol2003} N.\ Ueffing, L.\ Volkmann, 
   Restricted edge-connectivity and minimum edge-degree. 
   Ars Combin.\ {\bf 66} (2003), 193-203.
\bibitem{VolHon2017} L.\ Volkmann, Z-H Hong,
   Sufficient conditions for maximally edge-connected and super-edge-connected graphs.
   Commun.\ Comb.\ Optim.\ {\bf 2}(1) (2017), 35-41.
\bibitem{Whi1932} H.\ Whitney, 
   Congruent graphs and the connectivity of graphs. 
   Amer.\ J.\ Math.\ {\bf 54} (1932), 150-168.
\end{thebibliography}
\end{document}